# Развитие метода сечений Пуанкаре: визуализация трёхмерных сечений четырёхмерных потоков

А. Герега

*Научно-производственный центр, Одесса, Украина*

*E-mail: aherega@gmail.com*

Тема статьи – применение метода сечений Пуанкаре для визуальной классификации аттракторов в четырёхмерном фазовом пространстве; цель исследования – ввести в рассмотрение трёхмерные сечения Пуанкаре, и разработать алгоритм их использования для классификации четырёхмерных потоков по виду аттракторов, получаемых в таких сечениях.

В статье описано использование метода при изучении модели центробежного воздушного фильтра с обратными связями (ФОС): аппарата из последовательно соединенных криволинейных каналов (труб) постоянного сечения, в каждом из которых запылённый воздушный поток движется по дуге окружности, совершает поворот на 180°, и, благодаря взаимному поперечному смещению каналов, разделяется на две части. Это позволяет тяжёлым частицам, постоянно перемещаясь в каналы с большим радиусом кривизны, выйти из фильтра, а очищенному воздушному потоку попасть в атмосферу через торцевое отверстие в центре фильтра. В компьютерной модели ФОС рассматривается как открытая динамическая система с четырьмя взаимодействующими уровнями, соответствующими равновесным орбитам пыли. Модель описывается системой разностных уравнений; решение получено численно в виде аттракторов в четырёхмерном фазовом пространстве.

Моделирование динамики наполнения орбит частицами пыли показало, что в четырехмерном фазовом пространстве состояний каналов образуются аттракторы различных типов в зависимости от значений параметров, описывающих конструкцию и особенности функционирования фильтра. Статья содержит полученные в трёхмерных сечениях Пуанкаре следы $4D$-аттракторов, позволяющие составить представление об их топологии, а также провести их первичную классификацию. Кроме того, важно, что модель ФОС сформулирована как универсальная, что позволяет исследовать эволюцию систем (если они могут быть описаны итерационными уравнениями), состоящих из подсистем, взаимодействующих по произвольным законам.

*Ключевые слова*: динамическая система, хаос, классификация аттракторов, сечения Пуанкаре, топология, четырёхмерные потоки.

# Development of the Poincare cross-section method: Visualization the three-dimensional sections of four-dimensional flows

Alexander Herega

*Research and Production Center, Odessa, Ukraine*

*E-mail: aherega@gmail.com*

**Введение**

Классический метод сечения фазовой траектории в трёхмерном пространстве плоскостью, разработанный Анри Пуанкаре (Henri Poincare) [1], как известно, значительно упрощает исследование эволюции динамических систем. Причин несколько [1, 2]: во-первых, сечение Пуанкаре заменяет исследование эволюции трёхмерного потока с (квази)непрерывным временем изучением отображения на плоскости с дискретным временем; во-вторых, резко снижается количество данных, которые нужно обработать, так как «почти всеми точками траектории можно пренебречь» [1]. И третье, что особенно важно для исследуемой нами модели, – снижение на единицу размерности исследуемого объекта: при использовании метода мы переходим от потока в $R^n$ к отображению в $R^{n-1}$.

Развитие метода Пуанкаре рассматривается в статье на базе модели центробежного воздушного фильтра с обратными связями (ФОС). Эволюция состояния фильтра описывается в модели системой четырёх разностных уравнений, т.е. потоком в $R^4$. Понятны трудности описания такого потока: даже в трёхмерном пространстве траектории могут быть изображены лишь либо в перспективе, либо в проекциях, причём в обоих случаях бывает затруднительно получить представление об их конфигурации.

В предлагаемом в статье методе сечений для каждого четырёхмерного потока строится четыре трёхмерных сечения – вдоль каждой из координатных осей $R^4$. Это даёт возможность попытаться представить четырёхмерный поток по его трёхмерным проекциям, что конечно же непросто. (Кажется, здесь уместно вспомнить высказывание академика Л. Д. Ландау (Lev Landau), сделанное в 60-х годах прошлого века, что современная наука умеет исследовать явления, которые невозможно представить визуально).

**1. Динамический хаос в математической модели ФОС.
Предыстория**

Математическая модель центробежного воздушного фильтра, представленного как трёхуровневая открытая динамическая системы с взаимодействующими уровнями описана в [3-5].

Динамический фильтр с обратными связями (ФОС) представляет собой аппарат из последовательно соединенных криволинейных каналов (труб) постоянного сечения, в каждом из которых запылённый воздушный поток движется по дуге окружности, совершает поворот на 180°, и, благодаря взаимному поперечному смещению каналов, разделяется на две части: одна попадает в канал с меньшим радиусом кривизны, другая – с большим. Таким образом, конструкция ФОС позволяет тяжёлым частицам пыли, постоянно перемещаясь в каналы с большим радиусом кривизны, выйти из фильтра и быть уловленными, а очищенному воздушному потоку с остатками мелкой пыли выйти в атмосферу через торцевое отверстие в центре фильтра.

Очистка запыленного воздушного потока в ФОС осуществляется слоями пыли, циркулирующими по равновесным круговым орбитам. Именно в них частицы пыли вовлекаются в интенсивное взаимодействие, коагулируют и переходят на более «высокие» орбиты – в пылевые потоки других каналов и выводятся из фильтра. Таким образом, эффективность работы ФОС во многом определяется интенсивностью процесса коагуляции.

В компьютерной модели ФОС рассматривается как открытая динамическая система с тремя взаимодействующими уровнями, соответствующими равновесным орбитам в каналах. В модели исследуется характер эволюции системы, в зависимости от интенсивности входящего потока и конструктивных особенностей фильтра. Модель описывается системой уравнений

$$\Phi 1(x, y, z) = \begin{cases} x_{n+1} = x_n - k_{xy} p x_n^2 + k_{yx} q y_n^2 + x_{in} \\ y_{n+1} = y_n + k_{xy} p x_n^2 - \left(k_{yx} + k_{yz}\right) q y_n^2 + k_{zy} r z_n^2 \\ z_{n+1} = z_n + k_{yz} q y_n^2 - \left(k_{zy} + k_{out}\right) r z_n^2 \end{cases} \quad (1)$$

где $x, y, z$ – динамические переменные, определяющие количество частиц на уровнях, $k_{ij}$ – переходные коэффициенты, характеризующие статическое и динамическое взаимодействие уровней соответственно, $p, q, r$ – распределяющие коэффициенты, $x_{in}$ – количество частиц, входящих на первый уровень, причём, $\{k_{ij}\}$ и $\{p, q, r\} \in (0,1)$, $\{x, y, z\} \in R$, $x_{in} = const \in R^+$.

Наличие в системе двух групп коэффициентов $k_{ij}$ и $p, q, r$ объясняется её инженерным происхождением и имеет конкретную физическую интерпретацию: коэффициенты $k_{ij}$ – описывают относительное поперечное смещение последовательных каналов фильтра и задают долю потока, переходящего из одного канала в другой, а коэффициенты $p, q, r$ – описывают распределение частиц по ширине канала конструкции, так что переход между каналами определяется произведением коэффициентов обеих групп. Таким образом, каждое из итерационных уравнений системы является билинейным по $(p x_n) \cdot (k_{ij} x_n)$ [3-5].

Рассматриваемая система не интегрируется в общем виде, её решение может быть найдено лишь численно.

Стационарное решение системы уравнений, получаемое аналитически, имеет вид

$$x_{st} = \sqrt{\frac{x_{in}}{k_{xy} p}\left(1 + \frac{k_{yx}}{k_{yz}}\left(1 + \frac{k_{zy}}{k_{out}}\right)\right)} \quad y_{st} = \sqrt{\frac{x_{in}}{k_{yz} q}\left(1 + \frac{k_{zy}}{k_{out}}\right)} \quad z_{st} = \sqrt{\frac{x_{in}}{k_{out} r}} \quad (2)$$

С ростом $x_{in}$ возможны два варианта эволюции системы. В первом – возникает каскад бифуркаций удвоений периода, и реализуется сценарий Фейгенбаума (Feygenbaum) [6] перехода к хаосу. При этом в фазовом пространстве состояний равновесных орбит наблюдаются аттракторы, состоящие из $2^n$ точек. Во втором случае, после периодического режима возникает ситуация, аналогичная бифуркации Хопфа (Hopf) [7, 8], приводящей к возникновению квазипериодического режима и появлению в фазовом пространстве аттрактора в виде двух замкнутых линий. При реализации квазипериодического режима в модели может возникать синхронизация: тогда в системе наблюдается режим с зависящим от степени синхронизации периодом.

Дальнейшее увеличение $x_{in}$ приводит к бифуркации, в результате которой наступает хаотический режим, при этом в фазовом пространстве появляется странный аттрактор (рис. 1). Несколько иллюстраций: на рис. 2 показан аттрактор с фазовой траекторией, эволюционирующей от стационарного состояния (точка в центре) к аттрактору в виде кольца; на рис. 3 видна тонкая структура аттрактора. (Рисунки этого параграфа взяты из работы [5]).

Помимо рассмотренных сценариев развития хаоса могут наблюдаться различные их комбинации.

Дальнейшее возрастание $x_{in}$ в системе, уже находящейся в хаотическом режиме, приводит к возникновению перемежаемости; при этом в системе может наблюдаться мультиаттрактивность и гистерезис по параметру $x_{in}$.

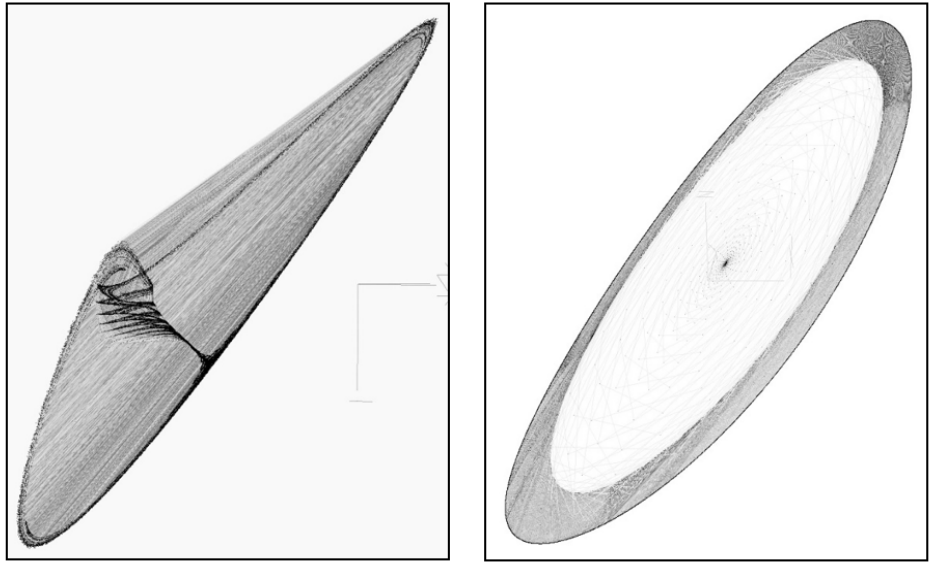

Рис 1. Странный аттрактор, возникающий в трёхмерном пространстве при параметрах $k_{xy}=0.5$, $k_{yx}=0.4$, $k_{yz}=0.3$, $k_{zy}=0.3$, $k_{out}=0.4$, $p=0.008$, $q=0.005$, $r=0.0057$, $x_{in}=39.65$.

Рис. 2. Аттрактор квазипериодического режима, возникающий в трёхмерном пространстве при параметрах ($k_{xy}=0.5$, $k_{yx}=0.1$, $k_{yz}=0.1$, $k_{zy}=0.4$, $k_{out}=0.5$, $p=0.05$, $q=0.02$, $r=0.01$, $x_{in}=30$.

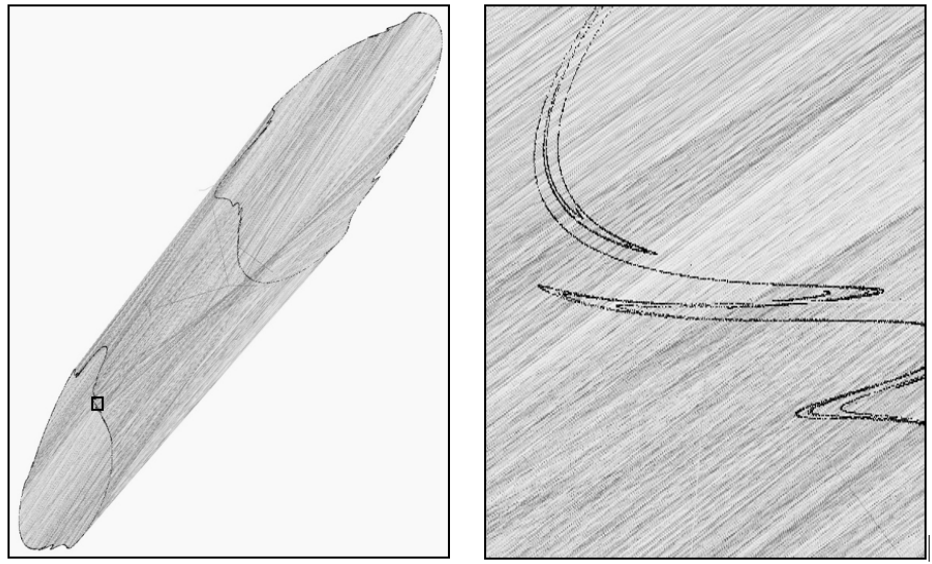

Рис. 3. Странный аттрактор, возникающий в трёхмерном пространстве при управляющих параметрах $k_{xy}=0.5$, $k_{yx}=0.1$, $k_{yz}=0.8$, $k_{zy}=0.2$, $k_{out}=0.4$, $p=0.05$, $q=0.05$, $r=0.045$, $x_{in}=11.6$; справа – увеличенный фрагмент аттрактора, выделенный прямоугольником.

## 2. Четырёхмерные потоки в модели центробежных фильтров

Рассмотрим модель центробежного фильтра с четырьмя равновесными орбитами. Модель описывается системой разностных уравнений

$$\Phi(x, y, z, w) = \begin{cases} x_{n+1} = x_n - k_{xy} p x_n^2 + k_{yx} q y_n^2 + x_{in} \\ y_{n+1} = y_n + k_{xy} p x_n^2 - (k_{yx} + k_{yz}) q y_n^2 + k_{zy} r z_n^2 \\ z_{n+1} = z_n + k_{yz} q y_n^2 - (k_{zy} + k_{zw}) r z_n^2 + k_{wz} s w_n^2 \\ w_{n+1} = w_n + k_{zw} r z_n^2 - (k_{wz} + k_{out}) s w_n^2 \end{cases}, \quad (3)$$

где $x$, $y$, $z$ и $w$ – динамические переменные, $k_{ij}$, $x_{in}$, $p$, $q$, $r$ и $s$ – распределительные коэффициенты, смысл которых описан выше.

Решение системы уравнений (3) получено в виде аттракторов различных типов в четырёхмерном фазовом пространстве состояний каналов (орбит), т.е. суммарной массы частиц в них. Тип аттракторов, разумеется, зависит от значений параметров, описывающих конструкцию и особенности функционирования фильтра.

Первичная классификация аттракторов, порождаемых четырёхмерными потоками, в компьютерной реализации модели проводится визуально по их 3$D$-следам (рис. 4-7), которые получены методом трёхмерных сечений.

При построении наиболее представительного, дающего максимальную информацию об аттракторе, сечения Пуанкаре в трёхмерном пространстве, удобно расположить систему координат так, чтобы секущая плоскость имела одну фиксированную координату. В четырёхмерном пространстве построение проводится аналогично: куб трёхмерного сечения, «перпендикулярный» одной из четырёх осей, также имеет одну фиксированную координату [9]. Существенное отличие в том, что в трёхмерном пространстве наиболее информативную плоскость несложно выбрать визуально, в четырёхмерном приходится «экспериментировать».

На рис. (4-7) показаны изображения странных аттракторов (полученных в четырёхмерном фазовом пространстве при определённых значениях параметров системы уравнений), в трёхмерных сечениях Пуанкаре вдоль различных осей, а также их фрактальные и корреляционные размерности. Относительная погрешность расчётов размерностей – несколько процентов.

## Заключение

В статье введены в рассмотрение трёхмерные сечения Пуанкаре, на использовании которых построен алгоритм анализа поведения динамических систем в четырёхмерном фазовом пространстве. Алгоритм позволяет получить представление о топологии их аттракторов по их виду в четырёх взаимно перпендикулярных трёхмерных сечениях Пуанкаре. Кроме того, визуализация следов таких аттракторов даёт возможность их систематизации по аналогии с визуальной классификацией трёхмерных аттракторов.

В дополнение к визуальной классификации в статье стандартным методом рассчитаны клеточные фрактальные размерности 4$D$-аттракторов, фрактальные и корреляционные размерности их 3$D$-следов.

Важно также, что исследуемая модель ФОС сформулирована как универсальная, и даёт возможность исследовать эволюцию любых систем (если они могут быть описаны итерационными уравнениями), состоящих из подсистем, взаимодействующих по произвольным законам.

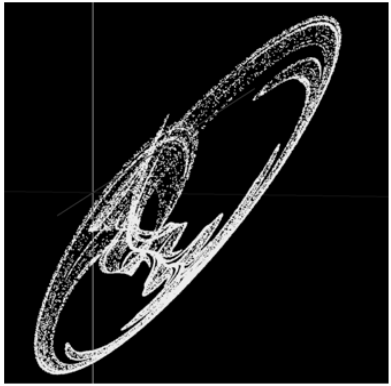
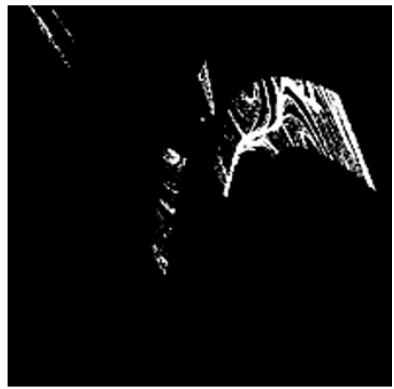
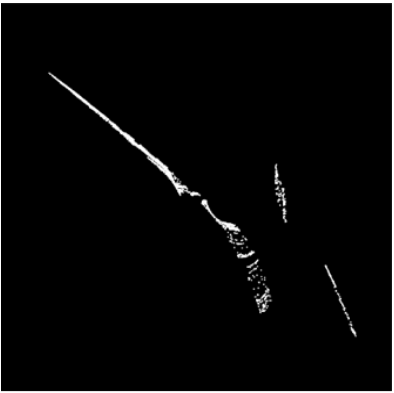
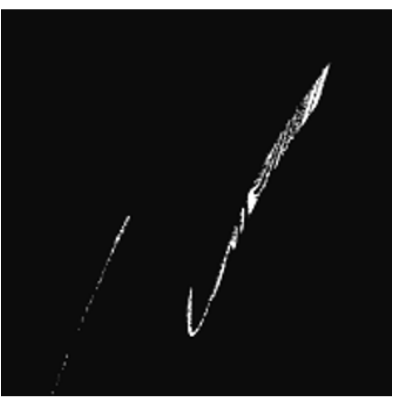

Рис. 4. Изображения странного аттрактора, полученного в четырёхмерном фазовом пространстве при значениях параметров системы уравнений: $k_{xy} = 0.65$, $k_{yx} = 0.25$, $k_{yz} = 0.65$, $k_{zy} = 0.25$, $k_{zw} = 0.000001$, $k_{wz} = 0.1$, $k_{out} = 0.4$, $p = 1$, $q = 1$, $r = 1$, $s = 0.1$, $x_{in} = 0.435$, в трёхмерных сечениях Пуанкаре вдоль пространственных осей $w$, $x$, $y$ и $z$ (сверху вниз).

Фрактальная размерность изображения в сечении перпендикулярном оси $w$ равна 2.56, корреляционная размерность 2.47. Фрактальная размерность 4*D*-аттрактора равна 3.58, корреляционная размерность 3.49.

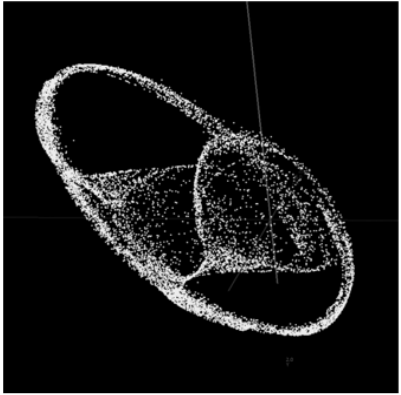

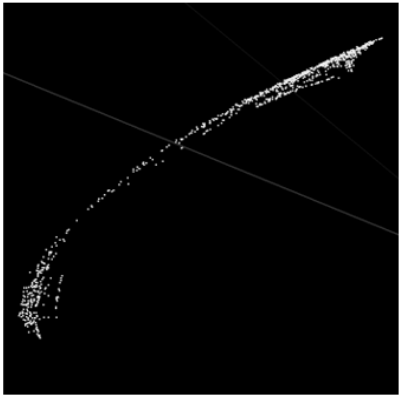

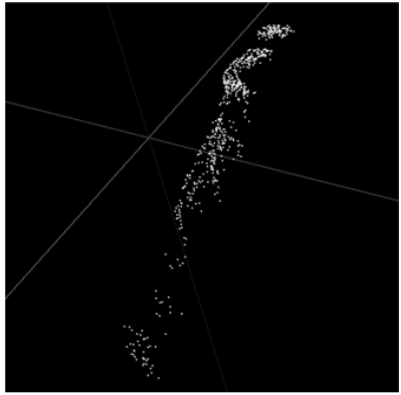

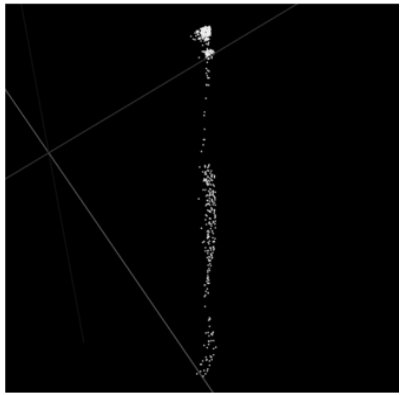

Рис. 5. Изображения странного аттрактора, полученного в четырёхмерном фазовом пространстве при значениях параметров системы уравнений $k_{xy} = 0.5$, $k_{yx} = 0.4$, $k_{yz} = 0.3$, $k_{zy} = 0.3$, $k_{zw} = 0.001$, $k_{wz} = 1$, $k_{out} = 0.4$, $p = 0.01$, $q = 1$, $r = 1$, $s = 1$, $x_{in} = 0.475$, в трёхмерных сечениях Пуанкаре вдоль пространственных осей $w$, $x$, $y$ и $z$.

Фрактальная размерность изображения в сечении перпендикулярном оси $w$ равна 2.56, корреляционная размерность 2.47. Фрактальная размерность 4D-аттрактора равна 3.66, корреляционная размерность 3.52.

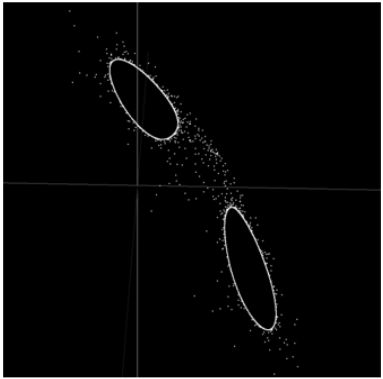
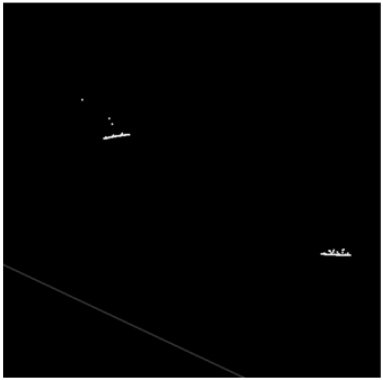
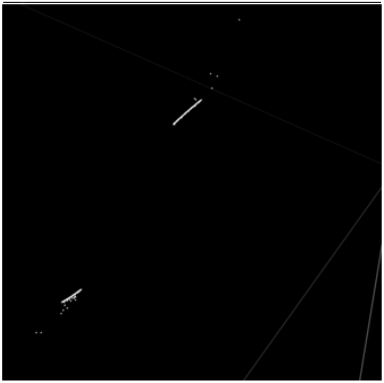
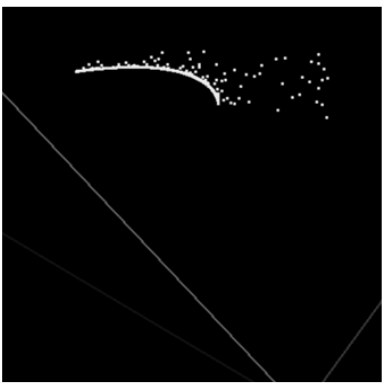

Рис. 6. Изображения странного аттрактора, полученного в четырёхмерном фазовом пространстве при значениях параметров системы уравнений $k_{xy} = 0.8$, $k_{yx} = 1.1$, $k_{yz} = 1.0$, $k_{zy} = 0.2$, $k_{zw} = 0.0001$, $k_{wz} = 0.1$, $k_{out} = 0.48$, $p = 0.95$, $q = 0.25$, $r = 0.2$, $s = 0.1$, $x_{in} = 0.33$, в трёхмерных сечениях Пуанкаре вдоль пространственных осей $w$, $x$, $y$ и $z$.

Фрактальная размерность изображения в сечении перпендикулярном оси $w$ равна 2.13, корреляционная размерность 2.04. Фрактальная размерность 4*D*-аттрактора равна 3.24, корреляционная размерность 3.19.

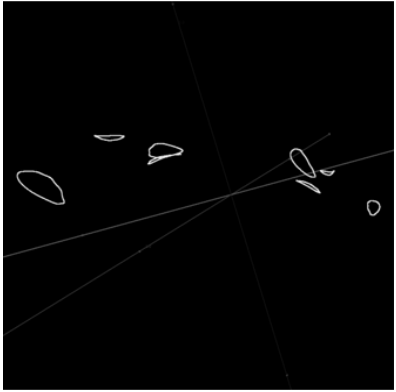
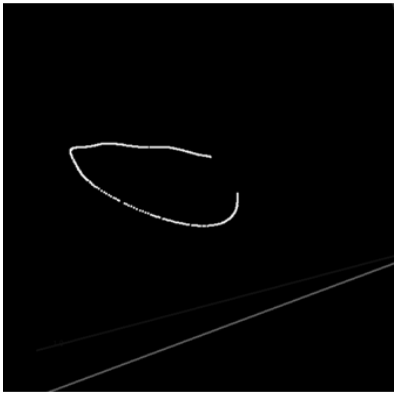
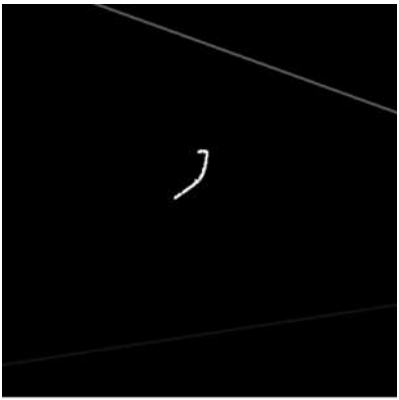
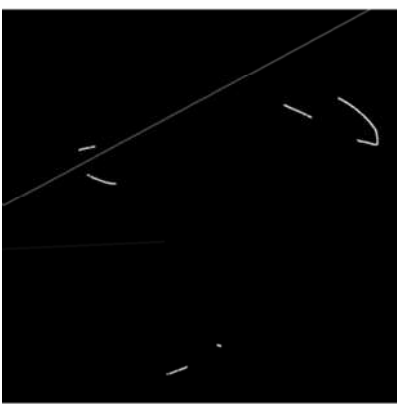

Рис. 7. Изображения квазипериодического аттрактора, полученного в четырёхмерном фазовом пространстве при значениях параметров системы уравнений $k_{xy} = 0.1$, $k_{yx} = 0.1$, $k_{yz} = 0.1$, $k_{zy} = 0.1$, $k_{zw} = 0.00001$, $k_{wz} = 0.1$, $k_{out} = 0.2$, $p = 1$, $q = 1$, $r = 1$, $s = 0.1$, $x_{in} = 1.536$, в трёхмерных сечениях Пуанкаре вдоль пространственных осей $w$, $x$, $y$ и $z$.

Фрактальная размерность изображения в сечении перпендикулярном оси $w$ равна 2.08, корреляционная размерность 2.03. Фрактальная размерность 4$D$-аттрактора равна 3.20, корреляционная размерность 3.14.